\documentclass[10pt]{article}
\usepackage{mathptmx}
\usepackage{color}

\input{epsf.tex}
\usepackage{latexsym,amsfonts,amssymb,epsfig,verbatim}
\usepackage{amsmath,amsthm,amssymb,latexsym,graphics,textcomp}
\usepackage{mathrsfs}
\usepackage{eucal}
\usepackage{graphicx}
\usepackage{url}

\input xy
\xyoption{all}

\theoremstyle{plain}
\newtheorem*{theo}{Theorem}
\newtheorem{theorem}{Theorem}
\newtheorem{lemma}[theorem]{Lemma}
\newtheorem{corollary}[theorem]{Corollary}
\newtheorem{example}[theorem]{Example}

\newtheorem{proposition}[theorem]{Proposition}

\theoremstyle{definition}

\newcommand{\PML}{\mathbb P \mathcal{ML}}

\newcommand{\AH}{\mathrm{AH}}

\newcommand{\Hom}{\mathrm{Hom}}

\newcommand{\Mod}{\mathrm{Mod}}

\newcommand{\PSL}{\mathrm{PSL}}

\newcommand{\SSS}{\mathbb S^3}
\newcommand{\BB}{\mathbb B}
\newcommand{\HH}{\mathbb H}

\newcommand{\CC}{\mathbb C}

\newcommand{\co}{\colon\thinspace}

\newcommand{\cC}{\mathcal{C}}
\newcommand{\CN}{\mathcal{N}}

\newcommand{\D}{\partial}
\newcommand{\Area}{\mathrm{area}}
\newcommand{\cusp}{\mathbf{P}}

\newcommand{\smallinfinity}{{\rotatebox{90}{\footnotesize $8$}}}

\title{\textbf{Geometric limits of knot complements, II: \\ Graphs determined by their complements}}
\author{Richard P. Kent IV and Juan Souto\thanks{The second author has been partially
supported by NSF grant DMS-0706878 and the Alfred P. Sloan Foundation.}}

\begin{document}

\date{April 5, 2009}

\maketitle

\begin{abstract}
\noindent We prove that there are compact submanifolds of the $3$--sphere whose interiors are not homeomorphic to any geometric limit of hyperbolic knot complements.
\end{abstract}

\bigskip
\noindent
A \textit{hyperbolic knot complement} is a complete hyperbolic $3$--manifold homeomorphic to the complement of a knot in $\SSS$. 
A complete hyperbolic $3$--manifold $M$ is a \textit{geometric limit} of hyperbolic knot complements if for every positive $\epsilon$, every compact submanifold $K$ in $M$ admits a $(1+\epsilon)$--bilipschitz embedding into a hyperbolic knot complement. 
Geometric limits of hyperbolic knot complements were studied by J. Purcell and the second author, who proved that every one--ended hyperbolic $3$--manifold with finitely generated fundamental group that embeds in $\SSS$ is a geometric limit of hyperbolic knot complements \cite[Theorem 1.1]{Jessica}.
It follows that every compression body is homeomorphic to a geometric limit of hyperbolic knot complements, and, in particular, there are geometric limits of hyperbolic knot complements with arbitrarily many ends.

The topology of such examples is not limited to that of compression bodies.

\begin{example}\label{example1}
There are compact hyperbolic $3$--manifolds with totally geodesic disconnected boundary  whose interior is homeomorphic to a geometric limit of hyperbolic knot complements.
\end{example}

It should be remarked that by \cite[Theorem 1.3]{Jessica}, a hyperbolic $3$--manifold with at least two convex cocompact ends is not a geometric limit of hyperbolic knot complements. 
The example shows that such a manifold may nonetheless be homeomorphic to one. 

These results lead Purcell and Souto to wonder if a compact submanifold of $\SSS$ whose interior admits a convex cocompact hyperbolic structure could fail to be homeomorphic to any geometric limit of knot complements \cite[Question 4]{Jessica}. We construct such a manifold here.

\begin{example}\label{main}
There is a compact submanifold of $\SSS$ that admits a hyperbolic metric with totally geodesic boundary whose interior is not homeomorphic to any geometric limit of knot complements.
\end{example}

Using elementary arguments about geometric limits, we obtain Example \ref{main} as a corollary of the following theorem.

\begin{theorem}\label{single-emb}
There is a compact oriented hyperbolic $3$--manifold $M$ with totally geodesic disconnected boundary that admits a unique orientation preserving embedding $M \to \SSS$ up to isotopy.
In fact, if $\Theta$ is a graph, then $\Theta$ admits an embedding into $\SSS$ so that the exterior of $\Theta$ has a unique oriented embedding into $\SSS$ up to isotopy.
\end{theorem}

\noindent Note that it follows from Fox's reembedding theorem \cite{Fox} that any compact submanifold of $\SSS$ that admits a single orientation preserving embedding in $\SSS$ is the exterior of an embedded graph. 

Using work of M. Lackenby \cite{Lackenby-attaching}, Theorem \ref{single-emb} is reduced to finding embeddings of $\Theta$ satisfying a certain geometric condition,  see Section \ref{shortobviousmeridians}. 
These embeddings are obtained in Sections \ref{reallyexcellent} and \ref{graphsshortmeridians} by a variation of arguments in \cite{Kent1}.
Theorem \ref{single-emb} and the examples are established in the final section.

\bigskip
\noindent
\textbf{Acknowledgment.} The authors thank Jessica Purcell for her interest and encouragement as this work developed, and 
Jeff Brock, 
Ken Bromberg,
and Dick Canary
 for nice conversation.

\section{Short obvious meridians and unique embeddings}\label{shortobviousmeridians}

A vertex of a graph is \textbf{extremal} if it has valence less than or equal to one.
In this section we consider finite embedded graphs $\Theta \subset \SSS$ without extremal vertices. 
A \textbf{branch vertex} is a vertex of valence at least three.
We say that $\Theta$ is \textbf{trivalent} if all of its branch vertices have valence three.

Every component of a regular neighborhood $\CN(\Theta)$ of an embedded graph $\Theta \subset \SSS$ is a handlebody. 
Given an edge $e$ of $\Theta$, we say that an essential simple closed curve $m\subset \D\CN(\Theta)$ is an \textbf{obvious meridian} associated to the edge $e$ if $m$ bounds a disk $D$ in $\CN(\Theta)$ that intersects $\Theta$ transversally in a single point in the interior of $e$. 
Since $\CN(\Theta)\setminus\Theta$ is homeomorphic to $\D\CN(\Theta)\times[0,1)$, every edge has an obvious meridian and any two obvious meridians for $e$ are isotopic in $\D\CN(\Theta)$. 
Note that if $\Theta$ is trivalent and $\Delta \subset \Theta$ is a connected component that is not a circle, then the obvious meridians form a pants decomposition of $\D \CN(\Delta)$.

Following Myers \cite{Myers}, a compact oriented $3$--manifold $M$ is \textbf{excellent} if it is irreducible atoroidal and acylindrical.
Equivalently, the complement of the torus boundary components of $M$ admits a complete hyperbolic metric with totally geodesic boundary, by Thurston's Geometrization Theorem for Haken Manifolds \cite{Otal-Haken}. 
By Mostow--Prasad rigidity \cite{mostow,prasad}, such a metric is unique up to isometry.
Given $\Theta$, we let $M_\Theta$ be the closure in $\SSS$ of the complement of $\CN(\Theta)$. 
With Myers, we say that $\Theta$ is \textbf{excellent} if $M_\Theta$ is. 

Suppose from now on that $M_\Theta$ is excellent and let $M_\Theta'$ be the complete hyperbolic manifold with totally geodesic boundary homeomorphic to the complement in $M_\Theta$ of the torus boundary components. 
Every torus  $T$ in $\D M_\Theta$ corresponds to a rank--two cusp of $M_\Theta'$. 
Let $\cusp(T)$ be the Margulis tube at this cusp with the property that the shortest essential curve in $\D\cusp(T)$ has length one. 
It was proved in \cite{BH} that the union $\cusp$ of these Margulis tubes over all torus components of $\D M_\Theta$ is a \textit{disjoint} union that does not meet $\D M_\Theta'$. 
We identify $M_\Theta$ with the complement of $\cusp$ in $M_\Theta'$. 
On each component of $\D M_\Theta$ of negative Euler characteristic this determines a hyperbolic metric, and on each torus a flat metric. 
So we may refer to the \textbf{length} $\ell(\gamma)$ of a curve $\gamma$ in $\D M_\Theta$ without risk of confusion.

Given $\epsilon$ positive, we say that an excellent graph $\Theta\subset\SSS$ has \textbf{$\epsilon$--short obvious meridians} if every obvious meridian $m$ on a torus $T$ in $\partial M_\Theta$ has length $\ell(m) < \epsilon\,\Area(T)$, and every other obvious meridian has $\ell(m)<\epsilon$.

We will refer to the subgraphs of $\Theta$ without extremal vertices as the \textbf{descendants} of $\Theta$.
Both $\Theta$ and $\emptyset$ are descendants of $\Theta$, and every descendant of a trivalent graph is trivalent. 
Our interest in excellent graphs, their descendants, and short obvious meridians is due to the following proposition.

\begin{proposition}\label{short-obvious}
For each $n$, there is an $\epsilon$ such that when $\Theta\subset\SSS$ is an excellent trivalent graph with $\vert\chi(\Theta)\vert\le n$ such that every nonempty descendant $\Delta \subset\Theta$ is excellent and has $\epsilon$--short obvious meridians, then every embedding $M_\Theta \to \SSS$ is the restriction of a diffeomorphism $\SSS\to\SSS$. 
In particular, $M_\Theta$ admits a unique orientation preserving embedding $M_\Theta \to \SSS$ up to isotopy.
\end{proposition}

We derive Proposition \ref{short-obvious} from the following result of Lackenby \cite{Lackenby-attaching}.

\begin{theo}[Lackenby]
Let $M$ be a compact excellent $3$--manifold and let $\cC$ be the set of essential simple closed curves $\gamma\subset\D M$ such that if $\gamma$ lies in a component $S$ of $\D M$ with negative Euler characteristic, then 
\[
\ell(\gamma) \leq \frac{4\pi}{(1-4/\chi(S))^{1/4}-(1-4/\chi(S))^{-1/4}}
\]
and if $\gamma$ lies in a torus, then $\ell(\gamma)\le 2\pi$.

Along each component $S$ of $\D M$ attach a $3$--manifold $H_S$ via a homeomorphism $S\to\D H_S$ to obtain a manifold $N$. 
Then either $N$ has infinite fundamental group or there is a component $S$ of $\D M$ containing a curve $\gamma \in \cC$ which is sent by $S\to\D H_S$ to a curve bounding a properly embedded disk in $H_S$.
\qed
\end{theo}
\noindent
In \cite{Lackenby-attaching}, this is only stated and proved when the $H_S$ are handlebodies,
but the proof applies without this restriction---see the proof of \cite[Theorem 1.3]{Jessica}. 
Proposition \ref{short-obvious} follows easily from Lackenby's theorem and the following observation.

\begin{lemma}\label{curves-curves}
For each $n$, there is an $\epsilon$ such that when $\Theta\subset\SSS$ is a trivalent excellent graph with $\vert\chi(\Theta)\vert\le n$ and $\epsilon$--short obvious meridians, then the collection $\cC$ in Lackenby's theorem is a collection of obvious meridians.
\end{lemma}
\begin{proof}
First consider a curve $\gamma$ in $\cC$ that lies in a component $S$ of $\D M_\Theta$ with $\chi(S)<0$.
 Note that $\vert\chi(S)\vert\le 2n$.  
Suppose that each obvious meridian in $S$ has length no more than $\epsilon$.
By the Collar Lemma,  each has a collar of width $\log\left(\coth\left(\epsilon/4\right)\right)$. 
 In particular, if $\epsilon$ is small enough that
\[
\log\left(\coth\left(\frac\epsilon 4\right)\right)>\frac{4\pi}{(1+4/n)^{1/4}-(1+4/n)^{-1/4}}
\]
then $\gamma$ must be isotopic off of the obvious meridians. 
Since the obvious meridians in $S$ form a pants decomposition, $\gamma$ is an obvious meridian itself.

Now consider a torus $T$ in $\D M_\Theta$ and recall that $T$ is endowed with a Euclidean metric of injectivity radius one. 
Suppose there is a curve $\gamma$ in $T$ with length $\ell(\gamma)\le 2\pi$ that is not isotopic to the obvious meridian $m$ in $T$. Then 
\[
1 \leq \Area(T)\le\ell(\gamma)\ell(m)\le 2\pi\ell(m),
\]
and this is not possible if $\Theta$ has $(1/4\pi)$--short obvious meridians. 
\end{proof}

\begin{proof}[Proof of Proposition \ref{short-obvious}]
Given $n$, let $\epsilon$ be the number provided by Lemma \ref{curves-curves} and suppose that $\Theta \subset \SSS$ is a trivalent graph such that all of its nonempty descendants are excellent and have $\epsilon$--short obvious meridians. 

Let $\varphi \co M_\Theta\to\SSS$ be an embedding. 
For every component $S$ of $\D M_\Theta$ let $H_S$ be the connected component of $\SSS\setminus\varphi(M_\Theta)$ facing $\varphi(S)$. 
So, gluing the manifolds $H_S$ to $N_{\Theta}$ via the boundary identifications induced by $\varphi$, we obtain the $3$--sphere. 
It follows from Lackenby's theorem and Lemma \ref{curves-curves} that there is an obvious meridian $m \subset \D M_\Theta$ whose image under $\varphi$ bounds a properly embedded disk in $H_S$. 
Let $e$ be the edge of $\Theta$ corresponding to $m$ and let $\Theta'$ be the largest descendant of $\Theta$ that does not contain $e$.

The manifold $M_{\Theta'}$ is homeomorphic to the manifold obtained by attaching a $2$--handle to $M_{\Theta}$ along $m$ and capping off boundary spheres. 
Since the image of $m$ under $\varphi$ bounds a disk in $H_S$, the embedding $\varphi$ extends to an embedding $\varphi'\co M_{\Theta'}\to\SSS$. If $\Theta'$ is empty, then $M_{\Theta'}=\SSS$ and we are done. 

Otherwise, observe that $\vert\chi(\Theta')\vert\le\vert\chi(\Theta)\vert$. 
By assumption, $\Theta'$ is excellent and has $\epsilon$--short meridians. 
The above argument yields a descendant $\Theta''$ of $\Theta'$ such that the embedding $\varphi'$ extends to an embedding $\varphi''\co N_{\Theta''}\to\SSS$. 
Repeating this process---no more often than the number of edges in $\Theta$---we obtain a diffeomorphism $\SSS\to\SSS$ which extends $\varphi$.

The final statement follows from the fact that the group of orientation preserving diffeomorphisms of $\SSS$ is connected.
\end{proof}

We now set about constructing graphs with arbitrarily short obvious meridians.

\section{Really excellent graphs}\label{reallyexcellent}

Let $X$ be a finite graph.
We let $\partial X$ denote the set of extremal vertices in $X$, and write $X^\circ = X - \partial X$.
Let $M$ a $3$--manifold.
An embedding $X \to M$ is \textbf{proper} if it induces a map of triples $(X, X^\circ\!, \D X) \to (M, M^\circ\!, \D M)$.
A finite graph $X$ properly embedded in $M$ is \textbf{excellent} if the exterior $M \setminus X$ is excellent.
We say that a proper embedding $X \to M$, or its image, is \textbf{nice} if $X$ is nonempty, has no isolated vertices, and no component of $\D (M \setminus X)$ is a sphere. 
Given a nice proper embedding $X\to M$,  a subgraph $Y$ is nice if the induced embedding is.

We will need the following mild generalization of Myers' theorem \cite{Myers}.
\begin{theorem}\label{general}
Let $M$ be an oriented $3$--manifold.
Let $X$ be a finite graph without isolated vertices. 
If there is a nice embedding $X \to M$,
then there is nice embedding $X \to M$ with the property that $M \setminus Y$ is excellent for any nice subgraph $Y$.
\end{theorem}

Theorem \ref{general} has the following corollary.

\begin{corollary}\label{totallynotbrunnian}  Let $X$ be a finite graph without isolated vertices properly embedded in a handlebody $H$ of positive genus.
Then there is an embedding $X \to H$ with the property that if $Y$ is any nonempty subgraph of $X$ without isolated vertices, then $H \setminus Y$ is excellent. \qed
\end{corollary}

We need some preliminary lemmata. 
Let $\BB$ be a $3$--ball. 
Recall that a tangle is \textbf{Brunnian} if removing \textit{any} strand results in the trivial tangle.

\begin{figure}
\begin{center}
\includegraphics{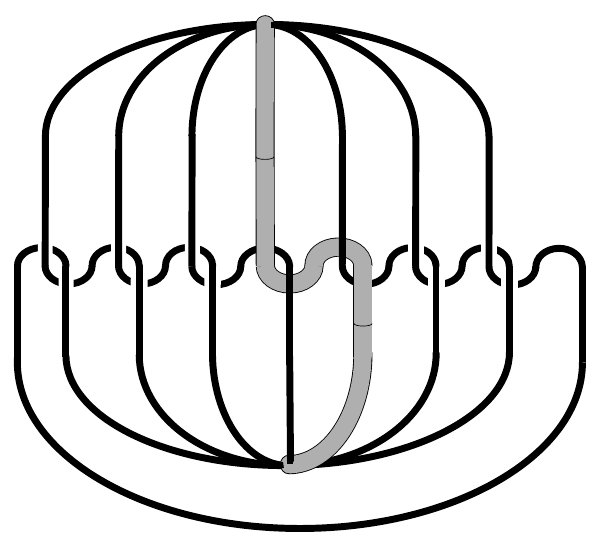}
\end{center}
\caption{An excellent Brunnian tangle on six strands.}\label{suzukitangle}
\end{figure}

\begin{lemma}\label{brunniantangle} For each $n$, there is an excellent Brunnian $n$--tangle in $\BB$.
\end{lemma}
\begin{proof} 
The exterior of Suzuki's Brunnian graph $\Theta_{n+1}$ on $n + 1$ edges, pictured in Figure \ref{suzuki}, admits a hyperbolic structure with totally geodesic boundary, see \cite{paoluzzizimmermann,ushijima}, and is thus excellent.
Moreover, the exterior of $\Theta_{n+1}$ is homeomorphic to the exterior of a Brunnian $n$--tangle in a ball---see Figure \ref{suzukitangle}.  
\begin{figure}
\begin{center}
\includegraphics{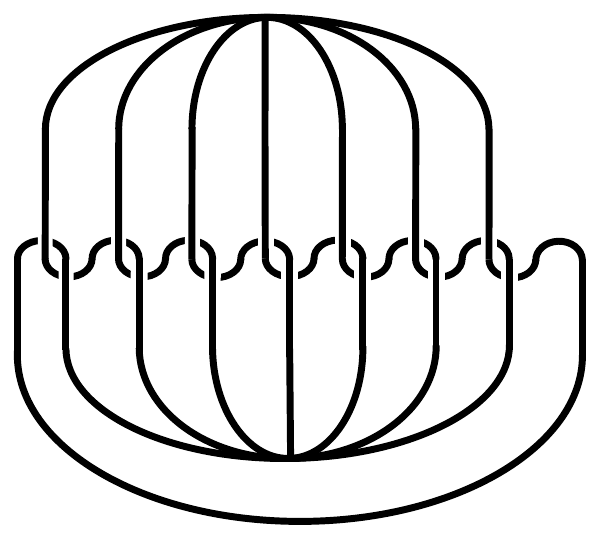}
\end{center}
\caption{Suzuki's Brunnian graph on seven edges.}\label{suzuki}
\end{figure}
Since the exterior of $\Theta_{n+1}$ is excellent, so is the exterior of the tangle.
\end{proof}

If $X$ is a graph, and $\Gamma$ a subgraph, we let  $X-\Gamma$ be the graph obtained by removing $\Gamma$ and taking the closure.
If $X$ and $Y$ are properly embedded graphs in a $3$--manifold $M$, we write $X \sim Y$ if they are ambiently isotopic.

\begin{lemma}\label{adaptation} Let $g \co X \to M$ be a nice embedding.
  Let $\Delta$ be a nice subgraph of $X$.
There is an embedding $f\co X \to M$  with $f(\Delta)$ excellent and such that for any edge $e$ of $\Delta$, 
\[
f(\Delta - e)  \sim g(\Delta - e)
\quad
\mathit{and} 
\quad
f(X - e) \sim  g(X - e).
\]
\end{lemma}
\begin{proof}
Let $e_1$, $e_2$, $\ldots$, $e_m$ be the edges of $\Delta$.
For each $i$, let $y_i$ be a point in the interior of $e_i$ and let $x_i = g(y_i)$.
Let $\alpha$ be an arc in $M$ whose left endpoint is $x_1$, whose right endpoint is $x_m$, and whose intersection with $g(X)$ is $\{x_1, \ldots, x_m\}$.

By Myers' theorem \cite{Myers}, the arc $\alpha$ is homotopic relative to $\{x_1, \ldots, x_m\}$ to an arc $\beta$ with $\beta \cap g(X) = \{x_1, \ldots, x_m\}$ such that $g(\Delta) \cup \beta$ is excellent.
 
Let $\CN(\beta)$ be a regular neighborhood of $\beta$ whose intersection with $g(X)$ is the standard trivial $m$--tangle in $\CN(\beta) \cong \BB$, each strand containing an $x_i$.
It is easy to see that $\D \CN(\beta) \setminus g(\Delta)$ is incompressible in $M \setminus (g(\Delta) \cup \CN(\beta))$.
Moreover, the exterior of $\CN(\beta) \setminus g(\Delta)$ in $M \setminus g(\Delta)$ is excellent, as it is homeomorphic to $M \setminus (g(\Delta) \cup \beta)$.

Now, in the obvious way, replace $\CN(\beta)  \setminus  g(\Delta) = \CN(\beta)  \setminus  g(X)$ with the exterior of an excellent Brunnian $m$--tangle given by Lemma \ref{brunniantangle}.
The resulting manifold is excellent, by Myers' gluing lemma, Lemma 2.1 of \cite{Myers}, and is homeomorphic to the exterior of an embedding $f\co \Delta \to M$.
Since the tangle is Brunnian, if $e$ is an edge of $\Delta$, we have 
\[
f(\Delta  - e) \sim  g(\Delta  - e)
\quad
\mathrm{and}
\quad
f(X  - e)  \sim  g(X  - e). \qedhere
\]
\end{proof}

\begin{lemma}\label{ordering} Let $X$ be a finite graph without isolated vertices and let $\mathfrak{S}(X)$ be the set of subgraphs of $X$ without isolated vertices. 
Then there is an ordering $\Gamma_1, \ldots, \Gamma_m$ of the elements of $\mathfrak{S}(X)$ such that $\Gamma_k \not\subset \Gamma_j$ when $k < j$.
\end{lemma}
\begin{proof} If $X$ is a single edge, there is nothing to do.

Let $E > 1$. 
Suppose that we have proven the lemma for all graphs whose number of edges is strictly less than $E$, and suppose that $X$ has $E$ edges.

Pick an edge $e$ of $X$.  
Partition $\mathfrak{S}(X)$ into two sets $\mathfrak{S}_+$ and $\mathfrak{S}_-$, the first consisting of those subgraphs that contain $e$, the latter of those that do not.
These sets are order isomorphic when equipped with the partial order induced by inclusion, as is easily seen by sending each subgraph in $\mathfrak{S}_-$ to the subgraph in $\mathfrak{S}_+$ obtained by adding $e$.

By induction, we may order $\mathfrak{S}_-$, and hence $\mathfrak{S}_+$, as desired.
Namely 
\[
\mathfrak{S}_+ = \{\Gamma_1, \ldots, \Gamma_\ell\}
\quad
\mathrm{and} 
\quad
 \mathfrak{S}_- = \{\Gamma_{\ell +1}, \ldots, \Gamma_{2\ell}\}.
 \]
As no element of $\mathfrak{S}_+$ is contained in any element of $\mathfrak{S}_-$, our desired ordering is $\Gamma_1, \ldots, \Gamma_{2\ell}$.
\end{proof}

\begin{proof}[Proof of Theorem \ref{general}]
Let $\Gamma_1, \ldots, \Gamma_N$ be the allowable \textit{removable} subgraphs of $X$---so each $X - \Gamma_i$ is nice.
By Lemma \ref{ordering}, we may assume $\Gamma_k \not\subset \Gamma_j$ when $k < j$, and we do so.

By Myers' theorem \cite{Myers}, 
there is an embedding $g_1 \co X \to M$ with image $X_1$ such that $X_1 - g_1(\Gamma_1)$ is excellent. 
To see this, note that there is an excellent embedding $g_0$ of $X - \Gamma_1$ in $M$.
We may then build an embedding of $\Gamma_1$ into $M - g_0(X - \Gamma_1)$ to obtain the desired embedding.

We want an embedding $g_N \co X \to M$ with excellent image $X_N$ so that $X_N - g_N(\Gamma_k)$ is excellent for all $k$.

Suppose that we have constructed an embedding $g_j \co X \to M$ with image $X_j$ so that for any $k \leq j$, the graph $X_j - g_j(\Gamma_k)$ is excellent.
We will construct an embedding $g_{j+1} \co X \to M$ with image $X_{j+1}$ so that for any $k \leq j+1$, the graph $X_{j+1} - g_{j+1}(\Gamma_k)$ is excellent.

Apply Lemma \ref{adaptation} to $g_j$ with $\Delta = X - \Gamma_{j+1}$ to obtain an embedding $f_{j} \co X \to M$ with image $Y_{j}$ such that $Y_{j} - f_{j}(\Gamma_{j+1})$ is excellent and deleting any edge $f_{j}(e)$ of $f_j(\Delta) =  Y_{j} - f_{j}(\Gamma_{j+1})$ from $Y_{j}$ yields a graph ambiently isotopic to $X_j -  g_j(e)$.

If $k < j + 1$, then $\Gamma_k \not\subset \Gamma_{j+1}$, and so there is an edge $e$ in $\Gamma_k$ that lies in $\Delta = X - \Gamma_{j+1}$.
So 
\[
Y_{j} - f_{j}(\Gamma_k) \, \sim \,  X_j - g_j(\Gamma_k)
\]
 when $k \leq j$, which is excellent by induction.

Now apply Lemma \ref{adaptation} to $f_j$ with $\Delta = X$
to obtain an embedding $g_{j+1}$ with excellent image $X_{j+1}$.
So deleting \textit{any} edge $g_{j+1}(e)$ from $X_{j+1}$ yields a graph ambiently isotopic to $Y_{j} - f_{j}(e)$.

So, if $k \leq j$, then $X_{j+1} - g_{j+1}(\Gamma_k)$ has the same exterior as $Y_{j} - f_{j}(\Gamma_k)$, which has the same exterior as $X_j - g_j(\Gamma_k)$, which is excellent by induction.

Finally note that $X_{j+1} - g_{j+1}(\Gamma_{j+1}) = Y_{j} - f_{j}(\Gamma_{j+1})$ is also excellent.

We are now done by induction.
\end{proof}


\section{Graphs with short obvious meridians}\label{graphsshortmeridians}

Let $F$ be a separating surface of genus at least two in a $3$--manifold $M$.
Let $\Theta$ be a graph intersecting $F$ transversely in a nonempty set of points such that the components $M_1$ and $M_2$ of $M \setminus(F \cup \Theta)$ are excellent.
Let $P$ be the union of tori in $\D M_\Theta$ and let $Q_j = P \cap M_j$.
We consider $(M,P)$ and the $(M_j,Q_j)$ as pared manifolds.
Let $S = F \setminus \Theta$.
So
$
M_\Theta = M_1 \cup_S M_2.
$
If $\varphi$ is a pure braid in $\Mod(S)$, let $M_\Theta^\varphi$ be the manifold obtained by gluing the $M_j$ together via $\varphi$.
Such an $M_\Theta^\varphi$ is the exterior of a different embedding of $\Theta$ into $M$ and admits a hyperbolic structure with totally geodesic boundary.
\begin{proposition}\label{everygraphpinch}
Let $\epsilon > 0$.
 If $\varphi$ is pseudo-Anosov, then for all sufficiently large $n$, the length of each component of $\partial S$ is less than $\epsilon$ in the metric on $M^n = M_\Theta^{\varphi^n}$ with totally geodesic boundary.
\end{proposition}
\begin{proof}
Equip $M^n$ with its hyperbolic metric with totally geodesic boundary, and let 
\[ 
\rho^n\co \pi_1(M^n)\to\PSL_2\CC
\] 
be the associated holonomy representation. 
Let 
$
\rho_j^n\co \pi_1(M_j)\to\PSL_2\CC
$
 be the representations induced by the inclusions $M_j \to M^n$. 
Let $\alpha$ be a simple closed essential and nonperipheral curve in $S$.

Since $M_2$ is excellent, 
the set 
\[
\AH(M_2) \subset \Hom(\pi_1(M_2),\PSL_2\CC)/\PSL_2\CC
\]
 of discrete faithful representations is compact, by \cite{thurstongeomI}.
In particular, there is an $L > 0$ such that, for all $n$,  the translation length $\ell(\rho_2^n(\alpha))$ of $\rho_2^n(\alpha)$ acting on $\HH^3$ is bounded above by $L$. 
It follows that, for all $n$, the geodesic representative of $\alpha$ in $M^n$ has length bounded by $L$. 
Viewing $S$ as a subsurface of $\D M_1$, we have
\begin{equation}\label{eq:length-bound}
\ell\left(\rho_{1}^n\left(\varphi^{-n}(\alpha)\right)\right) \leq L\ \ \hbox{for all}\ n.
\end{equation}

By \cite{thurstongeomIII}, every subsequence of the $\rho_1^n$ has a further subsequence that converges algebraically to the holonomy representation of a hyperbolic structure with parabolics at $Q_1$.
Let $\rho_1^\smallinfinity$ be such a limit.
In the space $\PML(S)$ of projective measured laminations on $S$, the sequence $\varphi^{-n}(\alpha)$ converges to the unstable lamination $\lambda^-$ of $\varphi$.
By \eqref{eq:length-bound} and continuity of Thurston's length function \cite[Corollary 7.3]{Brock-length-function}, the lamination $\lambda^-$ is not realized in $\rho_1^\smallinfinity$.
Thurston's compactness theorem for pleated surfaces \cite[Theorem 5.2.2]{notesonnotes} now implies that $\rho_1^\smallinfinity$ takes each component of $\D S$ to a parabolic element. 
Since this is true for any limit obtained by passing to a subsequence of the $\rho_1^n$, we have that,
 for all large $n$, the representation $\rho_1^n$ carries each component of $\D S$ to an element with small translation length.
 
It follows that if $\beta$ is a component of $\D S$ that does not lie in $Q_1$, its length in $M^n$ tends to zero.

Given a subsequence of the $\Gamma_1^{\, n} = \rho_1^n(\pi_1(M_1))$, we may always pass to a further subsequence that converges algebraically to some $\Gamma_1^{\, \smallinfinity}$ and geometrically to some $\widehat \Gamma_1$, see \cite[Corollary 9.18]{thurston} and \cite[Proposition 3.8]{jorgmarden2}.
The manifold $M_1^{\, \smallinfinity} = \HH^3/\Gamma^{\, \smallinfinity}$ covers $\widehat M_1 = \HH^3/\widehat \Gamma$.
By the above, the manifold $M_1^{\, \smallinfinity}$ has a degenerate relative end $E$ at $S$---see Section V of \cite{bonahon}.

Geometric convergence of a subsequence of the $\Gamma_1^{\, n}$ implies that we may pass to a further subsequence so that the $\Gamma^{\, n} = \rho^n(\pi_1(M^n))$ converge geometrically to a manifold $\widehat M$ covered by $\widehat M_1$, and hence $M_1^{\, \smallinfinity}$.
To see this, let $\gamma$ be a closed geodesic in $\widehat M_1$ corresponding to a nonperipheral embedded curve in $S$, and let $\gamma_n$ be the corresponding geodesics in the $M_1^n = \HH^3/\Gamma_1^{\, n}$.
We choose basepoints $x_n$ in the $\gamma_n$ to realize the geometric convergence $M_1^n \to \widehat M_1$, and let $y_n$ be the image of $x_n$ in $M^n$.
We claim that there is an $r > 0$ such that the injectivity radius of $M^n$ at $y_n$ is at least $r$.
Suppose that this is not the case.
Then for each $\epsilon > 0$, the image of $\gamma_n$ in $M^n$ intersects the $\epsilon$--thin part of $M^n$ for infinitely many $n$.
Now, the lengths of the $\gamma_n$ are uniformly bounded below.
Since the $M_1$ and $M_2$ are excellent, any embedded curve in $S$ is primitive in $\pi_1(M^n)$, and so $\gamma_n$ is primitive there.
It follows that $\gamma_n$ cannot lie entirely in an $\epsilon$--Margulis tube when $\epsilon$ is small compared to the length of $\gamma_n$.
By Brooks and Matelski's theorem \cite{brooksmatelski}, the distance between the boundary of the $\delta$--thin part and that of the $\epsilon$--thin part tends to infinity as $\epsilon$ tends to zero, and we are forced to conclude that the lengths of the $\gamma_n$ tend to infinity after passing to a subsequence.
But the lengths of the $\gamma_n$ are bounded above, as they tend to the length of $\gamma$.

Since the set of manifolds with injectivity radius $r >0$ at the basepoint is compact in the geometric topology \cite[Corollary 3.1.7]{notesonnotes}, we may pass to a subsequence so that the $(M^n, y_n)$ converge geometrically to a manifold $\widehat M$, which is covered by $\widehat M_1$.

By the Covering Theorem \cite{canarycovering}, the restriction of the covering map $M_1^{\, \smallinfinity} \to  \widehat M$ to $E$ is finite--to--one.
It follows that, in $\widehat M$, the cusps corresponding to the parabolic elements at $Q_i$ are rank--one cusps.

So, if $\beta$ lies in $\D Q_1 = \D Q_2$,
it has zero extremal length in a cusp cross section of any geometric limit $\widehat M$ of the $M^n$, and so its extremal length  tends to zero in the $M^n$.
\end{proof}


\begin{theorem}\label{shortmeridians} Let $\Theta$ be a nonempty graph with $\partial \Theta = \emptyset$ and let $\epsilon > 0$. 
There is an embedding of $\Theta$ into $\SSS$ such that if $\Delta$ is a nice subgraph of $\Theta$ , then $\SSS \setminus \Delta $ admits a hyperbolic structure with totally geodesic boundary and $\epsilon$--short obvious meridians.
\end{theorem}
\begin{proof}
Let $F$ be a genus--$2$ Heegaard surface cutting $\SSS$ into handlebodies $\mathbb{A}$ and $\mathbb{B}$.
Choose an embedding $\Theta \to \SSS$ such that the vertices of $\Theta$ miss $F$ and each edge of $\Theta$ intersects $F$ transversally in at least one point.
It follows that if $\Delta$ is a nonempty connected descendant of $\Theta$, then each edge of $\Delta$ intersects $F$.

The surface $F$ cuts $\Theta$ into two properly embedded graphs $X \subset \mathbb{A}$ and $Y \subset \BB$.
Let $X \to \mathbb{A}$ and $Y \to \BB$ be the embeddings given by Corollary \ref{totallynotbrunnian}, and call the images $X$ and $Y$ again.

Let $S = F - X = F - Y$, a punctured surface of genus two.
Given a pure braid $\zeta$ in $\Mod(S)$, we obtain a new embedding $\Theta \to \SSS$ by realizing $\zeta$ as a homeomorphism $h$ of $F$ fixing $\Theta \cap F$ pointwise, cutting $\SSS$ along $F$, and regluing via $h$.

Let $\Delta_1, \ldots, \Delta_n$ be the nice subgraphs of $\Theta$.
By reversing the order given by Lemma \ref{ordering}, we may assume that if $\Delta_j$ is a subgraph of $\Delta_i$, then $j < i$.
In particular, $\Delta_n = \Theta$.

For each $i$, let $S_i = F - \Delta_i$.
Each $S_i$ is a punctured surface of genus two, and so admits a Brunnian pseudo-Anosov braid $\varphi_i$, compare \cite{whittlesey}---a mapping class is \textbf{Brunnian} if it becomes the identity \textit{whenever} a puncture is filled.

If $\psi$ is a pure braid in $\Mod(S)$, we let $\kappa_j(\psi)$ denote the descendant of $\psi$ in $\Mod(S_j)$.
For each $j$, we choose a braid $\psi_j$ in $\Mod(S)$ with $\kappa_j(\psi_j) = \varphi_j$ and such that if $\Delta_j$ is not a subgraph of $\Delta_i$, then $\kappa_i(\psi_j)=1$.
This is possible as we may choose $\psi_j$ to be a braid descending to $\varphi_j$ such that filling any puncture of $S$ that survives in $S_j$ kills $\psi_j$.

If $a_1$, $\ldots$ , $a_{i-1}$ are integers, then, for all sufficiently large $N$, the pure braid
\[
\varphi_i^{N}  \kappa_i(\psi_{\, i-1}^{\, a_{i-1}}) \cdots \kappa_i(\psi_{\, 1}^{\, a_{1}})
\]
in $\Mod(S_i)$ is pseudo-Anosov.

Let 
\[
\zeta_\ell = \psi_{n}^{\, \ell_{n}} \cdots \, \psi_2^{\, \ell_2} \psi_1^{\, \ell_1},
\]
and let $\zeta_{\ell, i}= \kappa_i(\zeta_\ell)$.
By our ordering of the $\Delta_j$ and choice of the $\psi_j$,
\begin{align*}
\zeta_{\ell,i} &
= \kappa_i(\psi_{n}^{\, \ell_{n}}) \cdots \, \kappa_i(\psi_1^{\, \ell_1} ) \\
& = \varphi_i^{\ell_i}  \kappa_i(\psi_{\, i-1}^{\ell_{i-1}}) \cdots \kappa_i(\psi_{\, 1}^{\, \ell_{1}})
\end{align*}
Now, 
by Proposition \ref{everygraphpinch}, we may choose $\ell_n \gg \ell_{n-1} \gg \cdots  \gg \ell_1 \gg 1$ and $\mathbb{L} \gg \ell_n$ so that all of the gluing maps $\zeta_{\ell,i}^\mathbb{L} \co S_i \to S_i$ yield graphs $\Theta_i$ in $\SSS$ whose complements admit hyperbolic structures with totally geodesic boundary in which the obvious meridians all have length less than $\epsilon$.
So $\Theta_n$ is the image of the desired embedding.
\end{proof}

\section{Manifolds inaccessible to knot complements}\label{inaccesible}

\begin{proof}[Proof of Theorem \ref{single-emb}]
Let $\Theta$ be a graph and let $\Phi$ be the union of the tree components of $\Theta$.
A regular neighborhood of $\Phi$ is a disjoint union of balls, and it is easy to see that $\SSS \setminus \Theta$ has a unique oriented embedding into $\SSS$ if $\SSS \setminus(\Theta - \Phi)$ does.

So let $\Theta$ be a graph with no tree components. 
If $\Theta$ is empty, there is nothing to do.  
Assume that $\Theta$ is nonempty.
Every graph exterior in $\SSS$ is homeomorphic to the exterior of a $3$--valent graph, and so we assume that $\Theta$ is $3$--valent.
Since no component of $\Theta$ is a tree, its exterior is homeomorphic to the exterior of a $3$--valent graph without extremal vertices, and we assume that $\Theta$ has this property as well. 

Let $\epsilon$ be the constant provided by Theorem \ref{short-obvious} when $n=|\chi(\Theta)|$. 
By Theorem \ref{shortmeridians}, there is an embedding $\iota\co \Theta \to \SSS$ such that every nonempty descendant of $\iota(\Theta)$ is excellent and has $\epsilon$--short obvious meridians. 
Let $M=\SSS \setminus \CN(\iota(\Theta))$ be the exterior of a regular neighborhood of $\iota(\Theta)$ in $\SSS$. 
Since $\iota(\Theta)$ is excellent, $M$ admits a hyperbolic metric with totally geodesic boundary.
By Theorem \ref{short-obvious}, any two orientation preserving embeddings $M \to \SSS$ are isotopic. 
\end{proof}

We are now ready to produce the Examples from the introduction.

\begin{proof}[Example \ref{main}]
Let $N\subset\SSS$ be a manifold with disconnected boundary containing no tori or spheres provided by Theorem \ref{single-emb}.

Suppose there is a sequence of hyperbolic knot complements $M_i=\SSS\setminus K_i$ which converges geometrically to a complete hyperbolic manifold $M$ homeomorphic to the interior of $N$.
 Identify $N$ with a compact submanifold of $M$ in such a way that $M\setminus N$ has a product structure. 
Geometric convergence provides a sequence of better and better bilipschitz embeddings $\varphi_i\co N\to M_i$. 
The $M_i$ are knot complements and so sit naturally in the $3$--sphere. 
Composing the $\varphi_i$ with these natural embeddings $M_i \to \SSS$ we obtain embeddings $\psi_i \co N\to \SSS$. 
By Theorem \ref{single-emb}, we may precompose with a homeomorphism of $\SSS$ so that each $\psi_i$ is the identity map. 
It follows that, for all $i$, the knot $K_i$ is disjoint from $N$. 

Since $\D N$ is disconnected, $\SSS \setminus N$ has at least two components. 
Observe that by construction all these components are handlebodies. Passing to a subsequence we may assume that all of the $K_i$ lie in a single one of these, $U$ say. 
Letting $V$ be a component of $\SSS\setminus N$ different from $U$, we see that 
 if $m \subset \D N$ bounds a disk embedded in $V$ then the curves $\varphi_i(m)$ are homotopically trivial in the $M_i$ for all $i$. 

Now, the $\varphi_i$ are bilipschitz embeddings. 
So there is an $L > 1$ such that the curves $\varphi_i(m)$ have length at most $L$ for all $i$. 
A homotopically trivial curve of length $L$ in a hyperbolic $3$--manifold bounds a disk with diameter no more than $L$. 
It follows that the curve $m$ is homotopically trivial in the geometric limit $M$. 
But the inclusion $N\to M$ is a homotopy equivalence and $m$ is essential in $N$.
\end{proof}

\begin{proof}[Example \ref{example1}.] Consider the genus--$2$ Heegaard splitting $\SSS=H\cup H'$ and let $\Theta$ be an excellent graph in $H$. 
Let $\varphi\co \D H\to\D H$ be a pseudo-Anosov mapping class which extends to $H$ and whose attracting and repelling laminations lie in the Masur domain of $H'$, and consider the sequence of manifolds $N_n$ obtained by gluing $H\setminus\Theta$ and $H'$ via $\varphi^n$. 
Note that $N_n$ is homeomorphic to $\SSS\setminus\varphi^n(\Theta)$. 
As in \cite{Hossein}, it follows that for all sufficiently large $n$, the manifold $N_n$ is excellent and hence admits a hyperbolic metric with totally geodesic boundary.
Equip $N_n$ with this metric. 
It is an (advanced) exercise in Kleinian groups to prove that the sequence $N_n$ has a geometric limit $N_\smallinfinity$ homeomorphic to the interior of $H\setminus\Theta$. 
Furthermore, the manifold $N_n$ satisfies the conditions of \cite[Theorem 1.1]{Jessica} and is hence a geometric limit of knot complements. 
 Being a geometric limit of geometric limits of knot complements, $N_\smallinfinity$ is a geometric limit of knot complements itself. 
\end{proof}

\bibliographystyle{plain}
\bibliography{excellent}

\bigskip

\noindent Department of Mathematics, Brown University, Providence RI 02912 
\newline \noindent  \texttt{rkent@math.brown.edu}   

\bigskip
\noindent Department of Mathematics, University of Michigan, Ann Arbor MI 48109 
\newline \noindent  \texttt{jsouto@umich.edu}

\end{document}